\documentclass[twoside, a4paper, nofootinbib,11pt]{article}
\usepackage{amssymb}
\usepackage{IEEEtrantools}
\usepackage[mathscr]{eucal}
\usepackage[dvips]{graphicx}
\usepackage{amsmath}
\usepackage{amsthm}
\usepackage{pstricks}
\usepackage{caption}
\def \ni{\noindent}

\newcommand{\be}{\begin{equation}}
\newcommand{\ee}{\end{equation}}
\newcommand{\ben}{\begin{equation*}}
\newcommand{\een}{\end{equation*}}
\newcommand{\bes}{\begin{eqnarray}}
\newcommand{\ees}{\end{eqnarray}}
\newcommand{\besn}{\begin{IEEEeqnarray*}{rCl}}
\newcommand{\eesn}{\end{IEEEeqnarray*}}

\newcommand{\txt}{\textrm}

\newtheorem*{theorem*}{Theorem}

\newtheorem*{definition*}{Definition}
\newtheorem*{lemma*}{Lemma}
\newtheorem*{prop*}{Proposition}

\title{A Construction of Euclidean Invariant, Reflection Positive Measures on a Compactification of Distributions}
\author{T. Tlas}
\date{}

\begin{document}
\maketitle
\thispagestyle{empty}

\begin{abstract}
\ni A simple construction is given of a class of Euclidean invariant, reflection positive measures on a compactification of the space of distributions. An unusual feature is that the regularizations used are not reflection positive.
\end{abstract}

\vspace{0.5cm}

The goal of this paper is to give, under mild conditions, a very simple construction of a class of reflection positive, Euclidean invariant measures on a certain compactification of the space of distributions. The construction will work in any number of dimensions and for a wide class of local actions. We will restrict ourselves to a single real scalar field, but it will be clear that the arguments below can be easily extended to other situations. Roughly speaking, we will give a rigorous meaning to expressions of the following familiar form

\be
\label{eq:vague}
\frac{\int  e^{- S[\phi]} F[\phi]  d \phi}{\int  e^{- S[\phi]}   d \phi},
\ee

where $\phi$ is supposed to live in some space of ``functions'', $F[\phi]$ is a member of a useful class of functionals of $\phi$ (e.g. a trigonometric polynomial), and finally, $S[\phi] = \int (\nabla \phi)^2 + m^2 \phi^2 + \mathcal{L}[\phi] $, where $\mathcal{L}[\cdot]$ is a ``local'' functional of $\phi$ and its derivatives. \newline

Of course, it is well-known that the expression above as it stands is a mathematical fiction, since there is no useful way of giving meaning to the measure $d \phi$. Nonetheless, it is possible to proceed by combining the quadratic part of $S$ with $d \phi$ and work with the resulting Gaussian measure. Even then, the expression above is ill-defined due to the fact that $\mathcal{L}[\phi]$ is undefined on the support of the Gaussian measure. An enormous amount of work was expended to try to solve these difficulties, see \cite{jaffe, simon, strocci} and the numerous references therein. Roughly speaking, there are two kinds of problems one needs to deal with: the ones which appear at the short scales (ultraviolet) and those appearing at the long scales (infrared). One typically proceeds by regularizing the theory in some fashion, by putting in cut-offs and removing them in the end. We shall proceed in the same way. \newline

Let us give now the precise definitions. Fix $D \in \mathbb{N}$. This will be the number of dimensions in which the $\phi$'s will live in. It will be fixed throughout the paper. Let us first describe the infrared regularization. This is accomplished by essentially moving expression (\ref{eq:vague}) to a sphere. Thus, let $R >0$ and consider the sphere $\mathbb{S}_R$ in $\mathbb{R}^D \times \mathbb{R}$ given by the equation

\ben
x^2 + (y-R)^2 = R^2,
\een 

where $x \in \mathbb{R}^D$ and $y \in \mathbb{R}$. Let $\mathfrak{s}$ denote the stereographic projection from $\mathbb{S}_R- (0,2R)$ to $\mathbb{R}^D$ given by

\besn
\mathfrak{s}(x,y) = \frac{2R}{2R - y} x.
\eesn

Let $\Delta(D,l)$ stand for the dimension of the space of spherical harmonics of degree $l$.\footnote{$\Delta(D,l) = \frac{  (2l+ D - 1 )(l+ D - 2)! }{(D-1)!l! }  $, but we are not going to need the explicit expression in what follows.}  Let $\{Y_{l, m} \}_{m=0}^{\Delta(D,l)}$ be an orthonormal basis of the this space with respect to the $L^2$ product on $\mathbb{S}^D$, where the sphere is given the Hausdorff measure $\Omega_R$ induced from the Lebesgue measure on $\mathbb{R}^{D+1}$. If we denote by $\nabla^2$ the Laplace-Beltrami operator on the sphere, and recalling that $\nabla^2 Y_{l, m} = - l(l+D-1) Y_{l,m}$, we see that there is a unitary isometry between the Sobolev space\footnote{The most convenient definition of $H^{\alpha}$ for us is as the completion of the space of $C^\infty$ functions in the norm $|| (- \nabla^2 + 1)^\alpha f||_{L^2}$.} on the sphere of order $k$, $H^k$, and the set of all `sequences',\footnote{We are slightly abusing terminology here since these are labelled by two indices, but this should not cause any confusion.} $\{f_{l,m} : m=0, 1, \dots, \Delta(D,l); l=0, 1, \dots \}$ satisfying

\ben
\sum_{l=0, 1, \dots; m = 0, 1, \dots, \Delta(D,l)  } \bigg ( l(l+ D-1) +1 \bigg )^k |f_{l, m} |^2 < \infty.
\een 

It follows in turn that the dual of $H^k$, $H^{-k}$ is isometric with the space of sequences $\{  \phi_{l,m}: m=0, 1, \dots, \Delta(D,l);l = 0, 1, \dots \}$ satisfying

\ben
\sum_{l=0, 1, \dots; m = 0, 1, \dots, \Delta(D,l)  } \bigg ( l(l+ D-1) +1 \bigg )^{-k} |\phi_{l, m} |^2 < \infty.
\een 

Now, note that the expression

\be
\label{eq:cov}
B\Big ( \{f_{l,m} \} , \{ g_{l,m} \} \Big ) = \sum _{l=0, 1, \dots; m = 0, 1, \dots, \Delta(D,l)  } \frac{ \overline{f_{l,m}} g_{l,m} }{ l(l+ D-1) +1  }
\ee

defines a trace class bilinear form on any $H^k$ for a sufficiently large $k$. It thus follows by standard methods \cite{bogachev}, that (\ref{eq:cov}) is the covariance of a Gaussian measure $\mu$ supported on $H^{-k}$ for a sufficiently large $k$. At this point we select some such $k$ and will hold it fixed in what follows. Of course, the measure just described is a rigorous realization of the heuristic expression $e^{ - \int_{\mathbb{S}_R^D} \phi (\nabla^2 + 1) \phi  d\Omega_R  } d \phi$.\newline

We now move to the ultraviolet regularization. Let $h$ be a positive, smooth, compactly supported and rotationally invariant function on $\mathbb{R}^D$. For any $\Lambda >0$, let $\tilde{h}_\Lambda (x) = h\big (   \Lambda x  \big )$. Let

\ben
h_\Lambda (\theta) = \frac{ \tilde{h}_\Lambda \circ \mathfrak{s}(\theta)  }{ \int_{\mathbb{S}^D} \tilde{h}_\Lambda \circ \mathfrak{s}(\theta) d\Omega_R(\theta)  } .
\een

It is easy to see, e.g. by considering spherical harmonics expansions, that for any element $\phi \in H^{-k}$, we have that $\phi_\Lambda = h_\Lambda \ast \phi \in C^\infty$, where $\ast$ stands for the convolution on the sphere.\footnote{For example, one can define $f_1 \ast f_2 (\theta) = \int_{SO(D+1)} f_1(\theta) f_2(g \theta) dg$ where $dg$ is the Haar measure on the special orthogonal group  in $D+1$ dimensions, $SO(D+1)$.} Moreover, $\phi_\Lambda \to \phi$ (in the sense of distributions) as $\Lambda \to \infty$. Now, given any bounded measurable function $\mathcal{L}$ on $\mathbb{R}^N$, we have that the

\ben
\int_{\mathbb{S}_R} \mathcal{L} \bigg ( \phi_\Lambda (\theta) , \nabla^2 \phi_\Lambda, (\nabla^2)^2 \phi_\Lambda (\theta) , \dots, (\nabla^2)^N \phi_\Lambda (\theta)     \bigg ) d\Omega_R (\theta)
\een

is a well-defined, bounded function on $H^k$, which is invariant under the orthogonal group in $D+1$ dimensions, $O(D+1)$. \newline

Consider now the expression

\ben
\frac{ \int  F[\phi] e^{ - \int_{\mathbb{S}_R} \mathcal{L}\big  (\phi_\Lambda, \dots, (\nabla^2)^N \phi_{\Lambda} \big )  d\Omega_R  }   d\mu    }{   \int e^{-\int_{\mathbb{S}_R} \mathcal{L}\big  (\phi_\Lambda, \dots, (\nabla^2)^N \phi_{\Lambda} \big )  d\Omega_R}  d\mu   }.
\een

This is a well-defined version of (\ref{eq:vague}) for any bounded, measurable function $F$ on the support of $\mu$. One would like at this stage to send $R$ and $\Lambda$ to infinity. Therefore suppose $\{R_n\}_{n=1}^\infty$ and $\{\Lambda_n\}_{n=1}^\infty$ are two sequences of positive numbers with $R_n, \Lambda_n \to \infty$. Notice that if we replace $R$ and $\Lambda$ with $R_n$ and $\Lambda_n$, then the above expression gives a sequence of numbers whose absolute values are bounded from above by $|| F ||_{L^\infty}$. There are several straightforward ways to linearly extract a number from a bounded sequence. We're going to do it using a Banach limit $L$.\footnote{Another intuitively appealing procedure would be to use nonstandard analysis, by taking $R$ and $\Lambda$ unlimited, and then extracting the standard part of the limited expression above.} More precisely, let $L$ be an element of the dual of $l^\infty$ which is norm one, positive, shift invariant\footnote{We can relax the requirement of shift-invariance as nothing in the proofs below will depend on it.} and which coincides with the usual limit when it acts on a convergent $l^\infty$ element. The existence of such a functional $L$ is guaranteed by the Hahn-Banach theorem (see e.g. \cite{conway}). \newline

We are thus led to consider the following expression

\be
\label{eq:integral}
L \Bigg ( \frac{ \int  F[\phi] e^{ - \int_{\mathbb{S}_{R_n}} \mathcal{L}\big  (\phi_{\Lambda_n}, \dots, (\nabla^2)^N \phi_{\Lambda_n} \big )  d\Omega_{R_n}  }   d\mu    }{   \int e^{-\int_{\mathbb{S}_{R_n}} \mathcal{L}\big  (\phi_{\Lambda_n}, \dots, (\nabla^2)^N \phi_{\Lambda_n} \big )  d\Omega_{R_n}}  d\mu   } \Bigg ).
\ee

We will show momentarily that the expression above can be considered an integral with respect to a certain measure. However, since we're interested in measures which are Euclidean invariant and reflection positive, we need to restrict the class of functions $F$ one is willing to consider simply to make these concepts meaningful. We shall take it to be the class of cylindrical functions given in the following 

\begin{definition*}
Let $\mathcal{D}$ denote the space of smooth, compactly supported functions on $\mathbb{R}^D$ and $\mathcal{D}'$ be its dual, the space of distributions. A function $F$ on $\mathcal{D}'$ is said to be cylindrical if there is $m \in \mathbb{N}$, and there are $f_1, \dots, \allowbreak f_k \in  \mathcal{D}$ and a bounded continuous function $\tilde{F}$ on $\mathbb{R}^m$, such that $$F[T] =  \tilde{F}\Big(f_1(T), \dots, f_m(T) \Big ).$$ The set of all cylindrical functions will be denoted by $Cyl$.\end{definition*}

It is obvious that the set of cylindrical functions is a vector space. Moreover, for any cylindrical function $F[T] = \tilde{F}\Big(f_1(T), \dots, f_m(T) \Big )$, we have that

\ben
\breve{F} [\phi] =    \tilde{F}\Big(f_1 \circ \mathfrak{s} (\phi), \dots, f_m \circ \mathfrak{s}(\phi) \Big )
\een

is a bounded continuous function on $H^{-k}$. \newline

Before we state our main result there is one more issue we need to sort out. From renormalization group arguments, one in general would \emph{not} expect that the measure constructed from (\ref{eq:integral}) would be useful. This is because, as is familiar, one needs to adjust the bare parameters in the Lagrangian as the cutoffs are removed. This would correspond to making the function $\mathcal{L}$ above dependent on the cutoff, i.e. dependent on $n$. Moreover, one should also allow (e.g. by considering the case $D =2$) for $\mathcal{L}$ to be unbounded, at least in the limit \cite{simon}. Therefore, suppose that $\{\mathcal{L}_n \}_{n=1}^\infty$ is a sequence of bounded functions.\footnote{We are not assuming that the functions are uniformly bounded, nor that they are in fact functions of the same number of variables. In practice however, $\mathcal{L}_n$ is usually chosen to depend on a fixed, small number of variables. For example if one is trying to do $\phi^4$ theory, one choice is $$\mathcal{L}_n ( \phi_{\Lambda_n}, \nabla^2 \phi_{\Lambda_n}) = A_n f_n(\phi_{\Lambda_n} \nabla^2 \phi_{\Lambda_n}) + B_n g_n (\phi_{\Lambda_n}) + C_n h_n (\phi_{\Lambda_n}) ,    $$ where $f_n, g_n$ and $h_n$ are bounded functions which tend pointwise to $x, x^2$ and $x^4$ respectively as $n \to \infty$. Of course, $A_n, B_n$ and $C_n$ correspond to the usual bare constants (field strength, mass and coupling constant).} and consider the expression

\be
\label{eq:integral1}
L \Bigg ( \frac{ \int  F[\phi] e^{ - \int_{\mathbb{S}_{R_n}} \mathcal{L}_n\big  (\phi_{\Lambda_n}, \dots, (\nabla^2)^N \phi_{\Lambda_n} \big )  d\Omega_{R_n}  }   d\mu    }{   \int e^{-\int_{\mathbb{S}_{R_n}} \mathcal{L}_n\big  (\phi_{\Lambda_n}, \dots, (\nabla^2)^N \phi_{\Lambda_n} \big )  d\Omega_{R_n}}  d\mu   } \Bigg ) = I (F).
\ee

We can now state our main result in the following

\begin{theorem*}
There is a unique (up to homeomorphism) compactification $\mathring{\mathcal{D}'}$ of $\mathcal{D}'$ such that every cylindrical function $F$ on $\mathcal{D}'$ has a unique continuous extension $\mathring{F}$ to $\mathring{\mathcal{D}'}$. Also, there is a unique, rotationally invariant probability measure $\mu$ such that 
\ben
\int \mathring{F} d \mathring{\mu} = I(\breve{F}).
\een 

Moreover, one  can choose the sequences $\{R_n\}_{n=1}^\infty$ and $\{\Lambda_n\}_{n=1}^\infty$ such that $\mathring{\mu}$ is reflection positive. If, additionally, one has that $I \big ( || \phi ||_{H^{-k}}   \big ) < \infty$, then $\mathring{\mu}$ is invariant under translations as well.
\end{theorem*}

Note that usually \cite{jaffe}, Euclidean invariance and reflection positivity are defined for measures supported on $\mathcal{D}'$. However, in view of our choice of the class of functions which we're interested in integrating, we can use essentially the same definitions, which are:

\begin{itemize}
\item Euclidean invariance: For any element $E$ of the Euclidean group, and any cylindrical $F$, we have that

\ben
\int \mathring{F}_E d\mathring{\mu} = \int \mathring{F} d\mathring{\mu},
\een

where $F_E[\phi] = \tilde{F}\Big ( Ef_1(\phi), \dots, Ef_m(\phi)   \Big )$ and $Ef_j (x) = f_j(E^{-1}x)$.
\item Reflection positivity: Let $V^+$ stand for the subspace cylindrical functions such that the supports of $f_1, \dots, f_m $ are contained in the subset $(x_1, \dots, x_D) \in \mathbb{R}^D$ with $x_D > 0$. Let $\Theta: \mathbb{R}^D \to \mathbb{R}^D$ be the reflection in the $x_D$ coordinate. Then, for any $F \in V^+$, we have that

\ben
\int \Big ( \mathring{F} \mathring{F}_\Theta   \Big ) d \mathring{\mu} \geq 0.
\een
\end{itemize}

Before we proceed with the proof of the theorem, let us give a couple of remarks:

\begin{itemize}
\item Note that the condition $I\big ( || \phi ||_{H^{-k}} \big ) < \infty$ can be considered a very mild version of the analyticity axiom \cite{jaffe}, which in effect, would require that $I( e^{z \phi(f)} )$ is a holomorphic function of $z$. Intuitively, this corresponds to the constructed measure having an exponential fall off `at infinity' as opposed to the linear one required in the theorem. This, incidentally, would also guarantee that the support of $\mathring{\mu}$ is contained within $\mathcal{D}'$. This ties in with the remark below.
\item The theorem stated above gives a rather simple and a very general construction of measures on the \emph{compactification} of the space of distributions. It should be clear that since there is very little restriction on the form of the Lagrangians nor on their coupling constants one will be able to construct ``nontrivial'' measures (no matter how one chooses to define ``triviality''). Alas, this does not mean that all the usual difficulties of constructive quantum field theory are over. This is because, in this subject one is interested in constructing nontrivial measures on the space of distributions itself and not on its compactification. Thus, even with the theorem above one needs to do further work in order to show that the constructed measure is supported on the distributions (or at least its support has a nonempty intersection with them). This is similar to the issues in one of the proofs of Bochner-Minlos' theorem where one would like to show that the constructed measure has no support in the corona set (``at infinity'') \cite{herglotz}. Of course, it is not difficult, again using the vast freedom in the choice of the Lagrangians, to make sure that the constructed measure \textit{is} supported on the distributions (just send the ``coupling constants'' to zero sufficiently fast). However, unless one is extremely careful this will result in a Gaussian measure (even a Dirac one). It is precisely in trying to balance the two requirements, support on the distributions \textit{and} non-Gaussianity, that one needs renormalization group arguments. The theorem above does not address this point and only states that the result will be Euclidean invariant and reflection positive.
\end{itemize}

\begin{proof}

Recalling the definition of the topology on $\mathcal{D}'$, it follows at once that the set of cylindrical functions separates points from closed sets in $\mathcal{D}'$. This implies that one can imitate the standard arguments (see e.g. \cite{munkres}) used to show the existence and the properties of the Stone-C\v ech compactification, but with the algebra of all continuous functions being replaced with the algebra of cylindrical ones. This shows the existence and uniqueness of the compactification $\mathring{\mathcal{D}'}$ that we want, as well as the unique extension property for cylindrical functions.\newline

Let $A = \{ \mathring{F}: F \in Cyl \}$. It is obvious that $A$ is a subalgebra of continuous functions on $\mathring{\mathcal{D}'}$. Note that if $\{F_n \}_{n=1}^\infty$ is a sequence in $Cyl$ which converges uniformly on $\mathcal{D}'$, then $\{ \mathring{F}_n \}_{n=1}^\infty$ is also a uniformly convergent sequence on $\mathring{\mathcal{D}'}$. It is to check that $A$ vanishes nowhere and separates points on $\mathring{\mathcal{D}'}$. Thus, by Stone-Weierstrass, the uniform closure of $A$ coincides with $C(\mathring{\mathcal{D}'})$, the algebra of all continuous functions on $\mathring{\mathcal{D}'}$.\newline

Now, it is obvious that $\mathring{F} \to I(\breve{F})$ is a linear positive functional on $A$. Moreover, as was mentioned above, we trivially have that $\Big |  I(\breve{F}) \Big | \leq  || F ||_{L^\infty}$. It follows that the functional above extends uniquely to a linear positive functional on $C(\mathring{\mathcal{D}'})$. By Riesz-Markov, we have that there is a probability measure $\mathring{\mu}$ such that this functional coincides with the integral with respect to $\mathring{\mu}$. Uniqueness of $\mathring{\mu}$ follows from the fact that $A$ is dense in $C( \mathring{\mathcal{D}'})$.\footnote{This procedure of defining a measure on a space by going to the compactification was used in a different, simpler context in \cite{herglotz}. Also, a similar idea is utilized in the construction of the celebrated Ashtekar-Lewandowski measure, see e.g. \cite{thiemann, al}.}   \newline

Now, let $\breve{F}[\phi] = \tilde{F} \Big ( g_1 (\phi), \dots, g_m(\phi)  \Big )$. If $O \in O(D+1)$, let 

\ben
(\breve{F})_O[\phi] = \tilde{F} \Big ( O(g_1)(\phi), \dots, O(g_m)(\phi)  \Big ),
\een

where, as usual, $O(g)(\cdot) = g( O^{-1}(\cdot))$. \newline

Now using the fact that $\mu$ and $\int_{\mathbb{S}_{R_n}} \mathcal{L}_n \Big ( \phi_{\Lambda_n}, \dots \Big) d \Omega_{R_n}$ are $O(D+1)$ invariant, we have  

\ben
\int  (\breve{F})_O[\phi] e^{ - \int_{\mathbb{S}_{R_n}} \mathcal{L}_n\big  (\phi_{\Lambda_n}, \dots \big )  d\Omega_{R_n}  }   d\mu  = \int  \breve{F}[\phi] e^{ - \int_{\mathbb{S}_{R_n}} \mathcal{L}_n\big  (\phi_\Lambda, \dots \big )  d\Omega_{R_n}  }   d\mu  .
\een 

Now, if $O$ belongs to the $O(D)$ subgroup preserving the $y$ axis, it follows at once that $I(\breve{F}_O) = I ( (\breve{F})_O) =  I(\breve{F})$, and thus $\mathring{\mu}$ is rotationally invariant. What remains is to deal with translations and with reflection positivity. We shall handle reflection positivity first. The proof will, in effect, use the Markov property of the free quantum field  \cite{nelson}. The fact that a free quantum field on a Riemannian manifold with a reflection is reflection positive since it's Markovian was shown in \cite{dimock}.\footnote{The same fact was shown by different methods in \cite{dangelis, ritter} as well.} We shall, along the way, show essentially the same thing by a somewhat different route which is more convenient to our setting. \newline

We shall suppress the subscripts $n$ in $R$, $\Lambda$, $\mathcal{L}$ to reduce clutter. We will re-instate them later on when we'll deal with the limit.\newline

Now, let $\delta >0$.  

\ben
\mathbb{S}_R^{+\delta} = \mathbb{S}^R \cap \{   (x_1, \dots, x_{D-1}, x_D, y)\in \mathbb{R}^{D+1} :   x_D > \delta \},
\een

with $\mathbb{S}_R^{-\delta}$ having the same definition with the replacement $x_D < - \delta$. Also, let $\mathbb{S}_R^{0\delta} = \mathbb{S}_R - (\mathbb{S}_R^{+ \delta} \cup \mathbb{S}_R^{-\delta})$. Let $H^k_+$ be the closed subspace of $H^k$ which is the closure of $C^\infty$ functions supported in $\mathbb{S}_R^{+ \delta}$ with $H^k_-$ the analogous space for $\mathbb{S}_R^{-\delta}$, and let $H^k_0$ denote the orthogonal complement of $H^k_+ \oplus H^k_-$. It is trivial to see that the support of every element of $H^k_0$ is contained in $\mathbb{S}_R^{0\delta}$. Finally, let $P^\pm$ denote the orthogonal projections of $H^k$ onto $H^k_\pm$, with $P^0$ being the projection onto $H^k_0$. \newline

We extend now $\Theta$ to $\mathbb{R}^{D+1}$ in the obvious way, by keeping $y$ fixed, i.e. $\Theta (x_1, \dots, x_{D-1}, x_D, y) = (x_1, \dots, x_{D-1}, -x_D, y)$. It should be clear that $\Theta$ induces a unitary map, $f(\cdot) \to f(\Theta \cdot)$, from $H^k_+$ onto $H^k_-$. \newline

Now, observe that

\bes
\label{eq:decomp}
\nonumber
B \Big (    f, g \Big ) & = & B \Big ( (P^+ + P^0 + P^-) f, (P^+ + P^0 + P^-)g    \Big ) \\
\nonumber
& = & B\Big (P^+ f, P^+g \Big ) + B \Big (P^0 f, P^0 g \Big ) + B\Big (P^- f, P^- g \Big ) \\
& = & B^+ \Big(f, g\Big) + B^0\Big(f,g\Big) + B^-\Big(f,g\Big).
\ees

To see that there are no cross-terms above, consider e.g. $B\Big ( P^0 f, P^+ g \Big )$. From (\ref{eq:cov}), we see that it is equal to $\langle P^0 f, (- \nabla^2 + 1)^{-1} P^+ g \rangle_{L^2}$. We want to show that this expression vanishes. To do that, it is enough to show that $\langle P^0 f, (- \nabla^2 + 1)^{-1} h \rangle_{L^2} = 0$ for any $C^\infty$ function $h$ which is supported in $\mathbb{S}_R^{+\delta}$, as such functions are dense in $H^k_+$. Now, notice that there is a smooth function $\check{h}$ such that $h = (-\nabla^2 + 1) \check{h}$. Moreover, the support of $\check{h}$ is contained in $\mathbb{S}_R^{+\delta}$ (in fact in the support of $h$). Probably the easiest way to see this is to use the fact that 

\besn
0 & = & \int_{\mathbb{S}_R - \txt{supp}(h)} h \check{h} d\Omega_R = \int_{\mathbb{S}_R - \txt{supp}(h)}\Big ( |\nabla \check{h}|^2 + |\check{h}|^2 \Big ) d \Omega_R\\  
&\geq & \int_{\mathbb{S}_R - \txt{supp}(h)} |\check{h}|^2 d\Omega_R.
\eesn

Therefore, since the support of $P^0f$ is disjoint from that of $\check{h}$, we have that

\ben
\langle P^0 f, (- \nabla^2 + 1)^{-1} h \rangle_{L^2} = \langle P^0 f, \check{h} \rangle_{L^2} =0.
\een

The other cross terms are dealt with similarly. \newline

It is obvious that $B^+, B^0$, and $B^-$ in the decomposition (\ref{eq:decomp}) are symmetric, positive, and trace class. Moreover, they are supported on $H^k_+, H^k_0$, and $H^k_-$ respectively. It follows that there are three Gaussian measures $\mu^+, \mu^0$, and $\mu^-$, and a decomposition of the support of $\mu$ of the form $\phi = \phi^+ + \phi^0 + \phi^-$, with the corresponding supports of the Sobolev functions being in $\mathbb{S}_R^{+\delta}$, $\mathbb{S}_R^{0\delta}$ and $\mathbb{S}_R^{-\delta}$.   \newline

Now, let $F \in V^+$ with $F[T] = \tilde{F} \Big (f_1(T), \dots, f_m(T) \Big )$. Let $\delta = \frac{1}{\Lambda}$. Since the supports of $f_1, \dots, f_m$ are compact and are all contained in the half-space $\{(x_1, \dots, x_D) : x_D >0 \}$, it is clear that for all sufficiently large $R$'s and $\Lambda$'s one has that the supports of $f_1 \circ \mathfrak{s}, \dots, f_m \circ \mathfrak{s}$ are contained in $\mathbb{S}_R^{+ \delta}$.\footnote{If the supports of $f_1, \dots, f_m$ are contained in the subset $\{(x_1, \dots, x_D): A \leq x_D \leq B \}$, then the statement above is true provided one chooses e.g. $\delta < A$ and $R>B$. Eventually, we will have that $\delta \to 0$ and $R \to \infty$, so this will be satisfied.} Now let $\phi_\Lambda^\pm  = h_\Lambda \ast \phi^\pm$. It is clear from elementary geometry that one can choose $\alpha >0$ such that $\phi^+_\Lambda$ vanishes in $\mathbb{S}_R^{- (\alpha \delta)}$ with a symmetric statement for $\phi^-_\Lambda$.\footnote{If $R$ is very large, so that the stereographic sphere almost `coincides' with the plane, then it is clear that $\alpha \simeq 1$ should be sufficient. Thus e.g. $\alpha = 1000$ is more than enough for our purposes. Note that while eventually we will send $R$ and $\Lambda$ to $\infty$, $\alpha$ will be held fixed.  }

\besn
 & & \int  \breve{F}[\phi] \breve{F}_\Theta[\phi] e^{ - \int_{\mathbb{S}^{+(\alpha \delta)}_R \cup \mathbb{S}^{-(\alpha \delta)}_R } \mathcal{L}\big  (\phi_\Lambda, \dots \big )  d\Omega_R  }   d\mu[\phi]  =    \\
& & \int  \breve{F}[\phi^+] \breve{F}_\Theta[\phi^-] e^{ - \int_{\mathbb{S}^{+(\alpha \delta)}_R } \mathcal{L}\big  (\phi^+_\Lambda, \dots \big )  d\Omega_R  } e^{ - \int_{\mathbb{S}^{-(\alpha \delta)}_R } \mathcal{L}\big  (\phi^-_\Lambda, \dots \big )  d\Omega_R  }   d\mu^+[\phi^+] d\mu^-[\phi^-] \\
& & = \Bigg ( \int \breve{F}_\Theta[ \phi^-]  e^{ - \int_{\mathbb{S}^{-(\alpha \delta)}_R } \mathcal{L}\big  (\phi^-_\Lambda, \dots \big )  d\Omega_R  } d\mu^-[ \phi^-]    \Bigg ) \times \dots \\
& &\dots  \times  \Bigg ( \int \breve{F}[ \phi^+]  e^{ - \int_{\mathbb{S}^{+(\alpha \delta)}_R } \mathcal{L}\big  (\phi^+_\Lambda, \dots \big )  d\Omega_R  } d\mu^+[ \phi^+]    \Bigg ) = \\
 & = & \Bigg ( \int \breve{F}[ \phi^+]  e^{ - \int_{\mathbb{S}^{+(\alpha \delta)}_R } \mathcal{L}\big  (\phi^+_\Lambda, \dots \big )  d\Omega_R  } d\mu^+[ \phi^+]    \Bigg )^2, 
\eesn

The first equality above is a consequence of the decomposition of the measure $\mu$ just described, the fact that the integrand is independent of $\phi^0$, and that the

\besn
 \int_{\mathbb{S}^{+(\alpha \delta)}_R \cup \mathbb{S}^{-(\alpha \delta)}_R } \mathcal{L}\big  (\phi_\Lambda, \dots \big )  d\Omega_R & = & \\
 \int_{\mathbb{S}^{+(\alpha \delta)}_R} \mathcal{L} \big ( \phi^+_\Lambda, \dots \big ) d \Omega_R + \int_{\mathbb{S}^{-(\alpha \delta)}_R} \mathcal{L} \big ( \phi^-_\Lambda, \dots \big ) d \Omega_R & & 
\eesn

The final line is a consequence of the change of variables $\phi^- = \Theta \phi^+$ in the first term. To see that this is so, note that by a direct calculation, for any $f \in H^k_+$ we have that

\ben
\int e^{i \phi^+(f) } d\mu[\phi^+] = e^{- B^+(f , f)} = e^{- B^-(\Theta(f), \Theta(f)) }=  \int e^{i \phi^-( \Theta (f) )} d\mu[ \phi^-].
\een

From this, it follows by taking limits that for any element $F \in V^+$, we have that

\ben
\int \breve{F}[\phi^+ ] d\mu^+[\phi^+] = \int \breve{F}_\Theta[\phi^-] d\mu^-[\phi^-].
\een

If we now use that $\Theta \phi^+_\Lambda(x) = \phi^-_\Lambda(\Theta(x))$, then approximating $\mathcal{L}$ by a $C^\infty$ function, and then taking limits of Riemann sums and using dominated convergence we have what we want. \newline

Now, notice that (\ref{eq:integral1}) is invariant under $\mathcal{L} \to \mathcal{L} + \txt{constant}$. This means that without loss of generality, we can assume that $\mathcal{L} \geq 0$. Let $\mathfrak{M} = \sup_{x \in \mathbb{R}^N} \mathcal{L}(x)$. We then have that 

\besn
\bigg | e^{ - \int_{\mathbb{S}^{+(\alpha \delta)}_R \cup \mathbb{S}^{-(\alpha \delta)}_R } \mathcal{L}\big  (\phi_\Lambda, \dots \big )  d\Omega_R  }  - e^{ - \int_{\mathbb{S}_R } \mathcal{L}\big  (\phi_\Lambda, \dots \big )  d\Omega_R  }     \bigg | & \leq & \\
e^{ - \int_{\mathbb{S}_R } \mathcal{L}\big  (\phi_\Lambda, \dots \big )  d\Omega_R  }  \bigg | e^{\int_{\mathbb{S}^{0(\alpha \delta) } }\mathcal{L}\big  (\phi_\Lambda, \dots \big )  d\Omega_R } - 1   \bigg | & \lesssim & \\
e^{ - \int_{\mathbb{S}_R } \mathcal{L}\big  (\phi_\Lambda, \dots \big )  d\Omega_R  }  \int_{\mathbb{S}^{0(\alpha \delta) } }\mathcal{L}\big  (\phi_\Lambda, \dots \big )  d\Omega_R & \lesssim & \\
e^{ - \int_{\mathbb{S}_R } \mathcal{L}\big  (\phi_\Lambda, \dots \big )  d\Omega_R  } \Big (  \mathfrak{M} \delta R^{D} \Big ) =  e^{ - \int_{\mathbb{S}_R } \mathcal{L}\big  (\phi_\Lambda, \dots \big )  d\Omega_R  } \Big (  \mathfrak{M} \frac{ R^{D}}{\Lambda} \Big ) ,
\eesn

where the harmonic analysis notation $\lesssim$ above stands for ``less or equal than an irrelevant constant multiple of''. Putting back the subscripts $n$ we have thus shown that

\besn
\bigg | e^{ - \int_{\mathbb{S}^{+(\alpha \delta_n)}_R \cup \mathbb{S}^{-(\alpha \delta_n)}_{R_n} } \mathcal{L}_n \big  (\phi_{\Lambda_n}, \dots \big )  d\Omega_{R_n}  }  - e^{ - \int_{\mathbb{S}_{R_n} } \mathcal{L}_n\big  (\phi_{\Lambda_n}, \dots \big )  d\Omega_{R_n}  }     \bigg | & \lesssim & \\
e^{ - \int_{\mathbb{S}_{R_n} } \mathcal{L}_n\big  (\phi_\Lambda, \dots \big )  d\Omega_{R_n}  } \Big (  \mathfrak{M}_n \frac{ R_n^{D}}{\Lambda_n} \Big ) .
\eesn

Now, note that the above discussion works for any two sequences $\{R_n\}_{n=1}^\infty$ and $\{\Lambda \}_{n=1}^\infty$ as long as they go to $\infty$. In order to show reflection positivity, we need to choose our sequences\footnote{In other words, we're simply choosing the relative rate at which we're removing the ultraviolet and infrared cutoffs.} so that they satisfy 

\ben
\frac{\mathfrak{M}_n R^D_n}{\Lambda_n} \to 0.
\een

Putting everything together we thus have that

\bes
\label{eq:depend}
\nonumber
\Bigg | \frac{ \int  \breve{F}[\phi] \breve{F}_\Theta[\phi] \Big  ( e^{ - \int_{\mathbb{S}_{R_n}} \mathcal{L}\big  (\phi_{\Lambda_n}, \dots \big )  d\Omega_{R_n}  }   -  e^{ - \int_{\mathbb{S}^{+(\frac{\alpha}{\Lambda_n})}_{R_n} \cup \mathbb{S}^{-(\frac{\alpha }{\Lambda_n})}_{R_n} } \mathcal{L}\big  (\phi_{\Lambda_n}, \dots \big )  d\Omega_{R_n}  }    \Big )    d\mu         }{   \int e^{-\int_{\mathbb{S}_{R_n}} \mathcal{L}\big  (\phi_{\Lambda_n}, \dots \big )  d\Omega_{R_n}}  d\mu   } \Bigg | & &  \\
 \negmedspace{} \lesssim || F ||^2_{L^\infty} \mathfrak{M}_n \frac{R_n^D}{\Lambda_n}  \to  0. \qquad  \qquad \qquad \qquad \qquad  & & 
\ees

Therefore, we get that $I(\breve{F} \breve{F}_\Theta) \geq 0$ and thus we have reflection positivity. \newline

It remains to show invariance under translations. Suppose that $\mathfrak{t}$ is a translation by a vector $t$. Let $t^\perp$ stand for the subspace orthogonal to $t$ in $\mathbb{R}^D$. For every $R>0$, there is a unique $SO(D+1)$ rotation $O_{t, R}$ of $\mathbb{S}_R$ such that $\mathfrak{s}(O_{t, R}(0)) = t$ and $\mathfrak{s}(O_{t, R}( \mathfrak{s}^{-1}(t^\perp)))$ is orthogonal to $t$ at $t$.\footnote{Here, the reader should perhaps draw the case when $D = 2$.} \newline

Now, if $F[T] = \tilde{F} \Big ( f_1(T), \dots, f_m(T)    \Big )$ is cylindrical, such that $\tilde{F}$ is $C^1$ on $\mathbb{R}^m$, then

\besn
& & \Bigg | \int \bigg ( (\breve{F})_{O_{t, R_n}  } [\phi]  - \breve{F}_{\mathfrak{t} }  [\phi] \bigg ) e^{  - \int_{\mathbb{S}_{R_n} } \mathcal{L}_n\big ( \phi_{\Lambda_n}, \dots \big ) d\Omega_{R_n}     } d\mu \Bigg | \\  
 &\leq & || \tilde{F}'||_{L^\infty}  \int  \bigg ( \max_{j =1, \dots, m} \Big |\phi \Big ( (t f_j) \circ \mathfrak{s} - O_{t, R_n}(f_j)  \Big  ) \Big |     \bigg )  e^{  - \int_{\mathbb{S}_{R_n} } \mathcal{L}_n\big ( \phi_\Lambda, \dots \big ) d\Omega_{R_n}     }   d\mu  \\
 & & || \tilde{F}'||_{L^\infty} \max_{j =1, \dots, m} \Big | \Big | (t f_j) \circ \mathfrak{s} - O_{t, R_n}(f_j)  \Big  | \Big |_{H^k}  \times \dots \\
 & & \dots \times \int || \phi ||_{H^{-k}} e^{  - \int_{\mathbb{S}_{R_n} } \mathcal{L}\big ( \phi_{\Lambda_n}, \dots \big ) d\Omega_{R_n}     }   d\mu .
\eesn

Now, notice that $\big | \big |(t f_j) \circ \mathfrak{s} - O_{t, R_n}(f_j) \big | \big |_{H^k}$ goes to zero as $R_n \to \infty$. Then, if $I \big (   || \phi ||_{H^{-k}} \big ) < \infty$, and $\tilde{F}$ is $C^1$, we have that $I(\breve{F}_{\mathfrak{t}}) =  I(\breve{F})$. Since this equation holds on a dense subset of $Cyl$, it in fact holds everywhere, which concludes the proof.
\end{proof}

\textbf{Acknowledgments:} The author would like to thank J. Merhej for reading a preliminary version of this paper and for the numerous comments which greatly improved its readability.

\texttt{{\footnotesize Department of Mathematics, American University of Beirut, Beirut, Lebanon.}
}\\ \texttt{\footnotesize{Email address}} : \textbf{\footnotesize{tamer.tlas@aub.edu.lb}}

\end{document}